\documentclass[12pt]{article}
\usepackage[T2A]{fontenc}
\usepackage[utf8]{inputenc}
\usepackage[english]{babel}

\usepackage{graphicx}

\usepackage{amsmath}
\usepackage{amssymb}
\usepackage{amsthm}

\newtheorem{example}{Example}
\newtheorem{theorem}{Theorem}
\newtheorem{lemma}{Lemma}
\newtheorem{note}{Note}

\newtheorem{definition}{Definition}

\newtheorem{problem}{Problem}

\renewcommand{\det}{\mathop{\rm det}\nolimits}

\font\Bbb=msbm10 scaled 1200

\def\R{\hbox{\Bbb R}}

\def\I{\hbox{\Bbb I}}

\def\beq{\begin{equation}} \def\eeq{\end{equation}}
\font\gothic=eufm10
\def\goth#1{\hbox{\gothic #1}}

\begin{document}

\begin{center}
ANALOGUE OF PONTRYAGINS MAXIMUM PRINCIPLE \\
FOR MULTIPLE INTEGRALS MINIMIZATION PROBLEMS \\

M.I. ZELIKIN
\end{center}

\bigskip

\smallskip

 \begin{abstract} 

The theorem like Pontryagin's maximum principle for multiple integrals is proved. Unlike the usual maximum principle, the maximum should be taken not over all matrices, but only on matrices of rank one. Examples are given.

\end{abstract}

\bigskip

\section{INTRODUCTION}

\bigskip

\subsection{Statement of the problem and notations}

\bigskip

Let $\goth N$ be a domain on a smooth $n$-dimensional manifold. Let $\rho : \xi \to \goth N$ be a $\nu$-dimensional vector-bundle over the base $\goth N$, the fiber
of this bundle over the point $t \in \goth N$, that is the full preimage of the point $t$ of the mapping $\rho$, is the 
$\nu$-dimensional linear space. Local coordinates on 
$\goth N$ will be denoted by $t =(t^1,...t^n)$; that on fibers by $x = (x^1,...x^{\nu})$. The common convention on summation by repeating indices is used, therewith Latin indices, relating to coordinates of the base,run from 1 to $n$, while the Greek ones, relating to coordinates of fibers, run from 1 to $\nu$.  Multi-indices will be denoted by capital letters. The symbol $I$ means the full set from 1 to $n$. The symbol $\wedge$ will be sometimes omitted; it always is implicit in product of differentials. We denote by 
 $\I$ the identity matrix (its dimension is defined by the corresponding formula).

Let us consider the functional, its part related to the chart
 $V \subset \goth N$ is

\beq 
{\cal F} = \int_V f\left( t,x,\frac{Dx}{Dt}
\right) dt^I.  
\label{1}
 \eeq

Subsequent calculations will be made in coordinates of the chart 
$V$.

Let us denote by $J_1(\xi)$ the bundle of 1-jets over $\xi$,
and let

\beq
q^{\alpha}_i := \frac {\partial x^{\alpha}}{\partial t^i}=
g^{\alpha}_i(t,x)
\label{2}
\eeq

\noindent be a section of $J_1(\xi)$. The section can be considered as a field of slopes $\cal G$, i.e. as a distribution of $n$-dimensional planes in the space $\xi$.
We shall say that a manifold $M=\{x=\hat x(t)\}$ is imbedded into the field $\cal G$, if
$\frac{\partial \hat x^\alpha}{\partial t^i} = 
g^{\alpha}_i(t,x)$.

\begin{problem}

Given $y(\cdot) \in C^1(\partial \goth N)$, to find the strong minimum of the functional $(\ref{1})$ over $C^1$-manifolds $x(\cdot):\goth N \to \xi$ defined in $\goth N$ subject to the boundary conditions 
$x(t)|_{\partial \goth N}=y(t)$.  
\end{problem} 
 
\bigskip

\subsection{The second variation}

\bigskip

The natural necessary condition of minimum of the functional is the non-negativity of the second variation (\cite{B} - \cite{V}). Thorough investigation of the second variation for multiple integral is due to A.Clebsch. He studied Dirichlet functional

$$
\delta^2 {\cal F} =\int_V \left [ \frac {\partial^2 \hat f}
{\partial \left (\frac {\partial x^\alpha}{\partial t^i}\right)
\partial \left (\frac {\partial x^\beta}{\partial t^j}\right )}
\frac {\partial h^\alpha}{\partial t^i} \frac {\partial
h^\beta}{\partial t^j} + 2\frac {\partial^2 \hat f} {\partial
\left (\frac {\partial x^\alpha}{\partial t^i}\right) \partial
x_\beta} \frac {\partial h^\alpha}{\partial t^i} h^\beta + \frac
{\partial^2 \hat f}{\partial x^\alpha \partial x^\beta}h^\alpha
h^\beta \right ] dt^I.
$$

The symbol "hat" over a function (say, $\hat{f}$) means that into all its arguments one substitutes the extremal $\hat{x}(t)$.
By using ideas of many-dimensional Riccati techniques Clebsch presupposed the existence of solution to partial differential Riccati-type equation. Now we know that it is the many-dimensional counterpart of the condition of absence of the conjugate points. Using this solution he transforms variables and reduces the functional to its principal part, that is to the quadratic form from the converted first derivatives of desired functions. Hence, he almost proved the weak local minimum of the functional. It seems that Clebsch had in mind that for multiple integral it is valid the direct similarity of Legendre condition : The non-negativity of the second variation follows the non-negativity of the principal part of the quadratic form, defined on all $(n \times \nu)$-matrices
$$
q^\alpha_i = \frac {\partial x^\alpha}{\partial t^i}.
$$

In the work of J.Hadamard~$(\cite{H})$, which was published half a century after the work of Clebsch, it was shown that it is not truly. Hadamard clarified that the quadratic form of principal members is not obliged to be non-negative on the set of all matrices. He finds the following necessary condition of non-negativity of the second variation.

\begin{theorem}[Hadamard].

Let the functional

\beq \delta^2 {\cal F} =\int_V \left [a^{ij}_{\alpha \beta}(t)
\frac {\partial x^\alpha}{\partial t^i} \frac {\partial
x^\beta}{\partial t^j} + 2c^i_{\alpha \beta} (t) \frac {\partial
x^\alpha}{\partial t^i}x^\beta + b_{\alpha \beta}(t)x^\alpha
x^\beta \right ]dt^I
 \label{5}
 \eeq
\noindent be non-negative for $x(\cdot)$ that satisfies to boundary condition 
  
$$
x|_{\partial V} = 0.
$$

Then for all values of $t \in V$, the quadratic form 
$a^{ij}_{\alpha \beta}(t)q^\alpha _i q^\beta _j$ takes
non-negative values on all $(n \times \nu)$-matrices  of the form $q^\alpha _i = \xi^\alpha \eta _i$ (that is on all
matrices of rank one).
\end{theorem}

The theorem can be reformulated as follows:
The biquadratic form
$a^{ij}_{\alpha \beta}(t)\xi^\alpha \xi^\beta \eta_i\eta _j$ is non-negative for all values $t \in V$ and $\xi \in \R^\nu, \eta \in (\R^n)^*$.

The gap between necessary and sufficient conditions was essentially diminished by Van Hove~\cite{V}. At the end of 40-th years he proved that the natural amplification of the
Hadamard-Legendre condition:

\beq
 \frac {\partial^2 \hat f} {\partial
\left (\frac {\partial x^\alpha}{\partial t^i}\right) \partial
\left (\frac {\partial x^\beta}{\partial t^j}\right )} \xi^\alpha
\xi^\beta \eta^i \eta^j \ge \varepsilon |\xi|^2|\eta|^2
 \label{6}
 \eeq
\noindent is the locally sufficient condition of $C^1$-minimum. The expression "locally sufficient" means that the domain of integration is sufficiently small.  
 
The idea of the Van Hove's proof is as follows. Firstly, we freeze coefficients. That means that arguments 
$(t=t_0,x=x_0)$ in coefficients of the quadratic form standing under the sign of the integral~(\ref{5}) are fixed. 
(This does not affect on needed estimations because one can chose as the domain of integration an arbitrary small neighbourhood of the given point). Secondly. we apply to the functional obtained (that has constant coefficients) the Fourier transformation and then use the Parseval equality. The proof provides an explanation of the fact -- why the Hadamard-Legendre condition uses only matrices of the rank one. Namely. the operator of differentiation transits under the Fourier transform into the operator of multiplication by independent variable. If the Fourier-image of the function $x^\alpha (t)$ is $\xi^\alpha (\eta)$, then the Fourier-image of its derivative
$\frac {\partial x^\alpha}{\partial t^i}$ will be
$\xi^\alpha \eta_i$. As a result, in the integrand of the Parseval equality arises biquadratic form and the inequality
(\ref{6}) assures the positivity of the variation.

The condition~(\ref{6}) relates to the notion of hyperbolicity of non-linear Euler systems for non-stationary processes. It provides the correctness of Cauchy problem for systems with the same coefficients under higher derivatives 
~\cite{Fr}, \cite{Ho}.

We denote by $\cal U$ the set of matrices of the rank 1.

\smallskip

Considering candidates for extension of functional like that of~(\ref{1}), it will be interesting to use the space $BV$ -- functions of bounded variation on the space $V \subset \R^n$.
Let $f \in L^1(V)$. We will consider $\int |Df|dx$ in the sense of distributions as $\sup \{ \int_V f div\, g\}dx$, where $g \in C_0^1(V,\R^n),\, \parallel g\parallel \le 1\}$. Such extensions
was used in the theory of minimal surfaces and in plasticity theory (see, for example, $\cite{G}, \cite{Gi}, \cite{MM}$).

\begin{definition}

The function $f$ is the function of bounded variation if
$$
\int_V |Df|dx < \infty.
$$
\end{definition}

The following claim is evident:

Let $V \in \R^n$ is open and $f_i \in BV(V)$ converges in
$L^1$ to a function $f$. Then  
$$
\int_V |Df|dx \le \lim_{j \to \infty} inf \int_V |Df_j|dx.
$$

The space $BV(V)$ is the Banach-space with the norm
$$
\parallel f\parallel_BV = \parallel f\parallel_{L^1} 
+ \int_V |Df|dx. 
$$

Let us note that $BV$ contains characteristic functions
$\varphi({\cal E})$ of open bounded sets ${\cal E} \in \R^n$
with smooth boundaries. Indeed, denote by
$\mu_{n-1}({\cal E})$ the $(n-1)$-dimensional measure of the boundary of the set ${\cal E}$. Then
$$
\int_{\R^n}|D\varphi ({\cal E})|dx =  
\mu_{n-1}(\partial {\cal E} \cap V)
$$

Hence, $BV(V)$ contains discontinuous functions and thus using of this class of functions would need additional restrictions on  growth of integrand $f$. Using the discontinuous variations would give much more simple proof of the maximum principle. But we prefer to restrict ourselves only by the assumption of the strong minimum that is achieved on a separate trajectory.

To prove the main theorem of the work it will be sufficient
to consider the extension of {\cal F} that is obtained by adding to
$C^1(V)$ the following rather narrow space of piece-wise linear functions $\goth X$. Let the domain $\Omega$ in the space $(t)$ be bounded by a polyhedron $\partial \Omega$. Let $\{\cup \Xi_i\}$ be the simplicial division of 
$\partial \Omega$. The space $\goth X$ consists of continuous
piece-wise linear functions ${\goth y}:\partial \Omega \to \R^n$ with the support in $\Xi_i$ and with values in one-dimensional space 
$\goth Y \subset \R^n$. The value of the integral {\cal F}
equals to sum of its values over all simplexes. It is clear that such functions can be approximated by smooth ones in
$C^1$-metric. It is sufficient to smooth out angles between any two adjacent simplexes using cylindrical surface with the generator parallel to its intersection, and then smoothly to seam these cylinders.

In the construction of Weyl $(\cite{W})$, it is used the trace of the product of matrices as a scalar product. So,
the part of Pontryagin's function $H$ will play

$H=-f+q^i_{\alpha}\frac{\partial x^{\alpha}}{\partial t^i}$. 

\begin{note}

Other forms of the scalar products (scalar products in spaces of external degrees) give another forms of Weierstrass functions that was introduced and used in the work (\cite{MZ}).
\end{note}

\bigskip

The Weyl's canonical system for multiple integrals has the form
\beq
\left \{
\begin{array}{rcl}
\frac{\partial x^{\alpha}}{\partial t^i}=& 
(\frac{\partial H}{\partial q}):=& 
\phi_i^{\alpha} \\
\frac{\partial q^i_{\alpha}}{\partial t^i}=&
-(\frac{\partial H}{\partial x})=& 
\frac{\partial f}{\partial x^{\alpha}}-
q^j_{\beta} \frac{\partial \phi_j^{\beta}}{\partial x^{\alpha}} 
\end{array}
\right.
\label{10}
\eeq

Take the variation of the first group of $(\ref{10})$ by using the variation $(\delta x)(t)$. We have

\beq
\frac{\partial (\delta x)^{\alpha}}{\partial t^i}=
\frac{\partial \phi_i^{\alpha}}{\partial x^{\beta}}(\delta x)^{\beta}
\label{11}
\eeq

The system $(\ref{11})$ is conjugated to the homogeneous part of the second group of equations $(\ref{10})$. As in the case of one-dimensional integral, the impulses $q^i_{\alpha}$  appear to be the infinitesimal tangent planes to the level surface of solutions. Thus the scalar product of $q^i_{\alpha}$ -- the solution to the equation
$$
\frac{\partial q^i_{\alpha}}{\partial t^i}= -
q^j_{\beta} \frac{\partial \phi_j^{\beta}}{\partial x^{\alpha}} 
$$
\noindent and the solution of the variational equation $(\ref{11})$  remains constant. Indeed,

\beq
\frac{\partial }{\partial t^i}[q^i_{\alpha}(\delta x)^{\alpha}]=
q^i_{\alpha}\frac{\partial \phi^{\alpha}_i}{\partial x^{\beta}}
(\delta x)^{\beta}-
q^i_{\alpha}\frac{\partial \phi^{\alpha}_i}{\partial x^{\beta}}
(\delta x)^{\beta}=0
\label{12}
\eeq

The immediate sequence of $(\ref{12})$ is the following lemma:

\begin{lemma}[On conjugate systems]

Let $\hat x(t)$ be an optimal solution and $(\hat x,\hat q)$ be the solution of Weyl's canonical system $(\ref{10})$.

Then,
  
$\frac{\partial }{\partial t^i}[\hat {q}_{\alpha}^i
(\delta x)^{\alpha}]=
\frac{\partial \hat {f}}{\partial x^{\alpha}}(\delta x)^{\alpha}$.
\end{lemma}

\smallskip

\subsection{The statement of the main theorem}

\smallskip

\begin{theorem}[The maximum principle]

Suppose that the function $f$ is smooth.
Suppose that $\hat x(\cdot)$ gives the strong minimum to the 
problem 1 with the functional.

Then there exists a solution to conjugate system of equations 
$(\ref{10})$ such that the Pontryagin's function $H$ attains its maximal value on the set of slopes $\frac{Dx}{Dt}$ defined by the rank one matrices -- $\cal U$:
\beq
\max_{(\xi^{\alpha} \eta^i) \in  {\cal U}}
\left[ -f\left( t,\hat x, \frac{D\hat{x}}{Dt}+(\xi^{\alpha} \eta_i)
\right) + \frac{\partial f}{\partial (\frac{Dx}{Dt})} 
(\xi^{\alpha} \eta_i)\right] = -\hat f +
\frac{\partial \hat f}{\partial (\frac{Dx}{Dt})}
\left( \frac{D\hat x}{Dt}\right) .
\label{13}
\eeq
\end{theorem}

\bigskip

\bigskip

\section{PROOF OF THE THEOREM}

\bigskip

\begin{figure}
 \centering
 \includegraphics[width=0.7\textwidth]{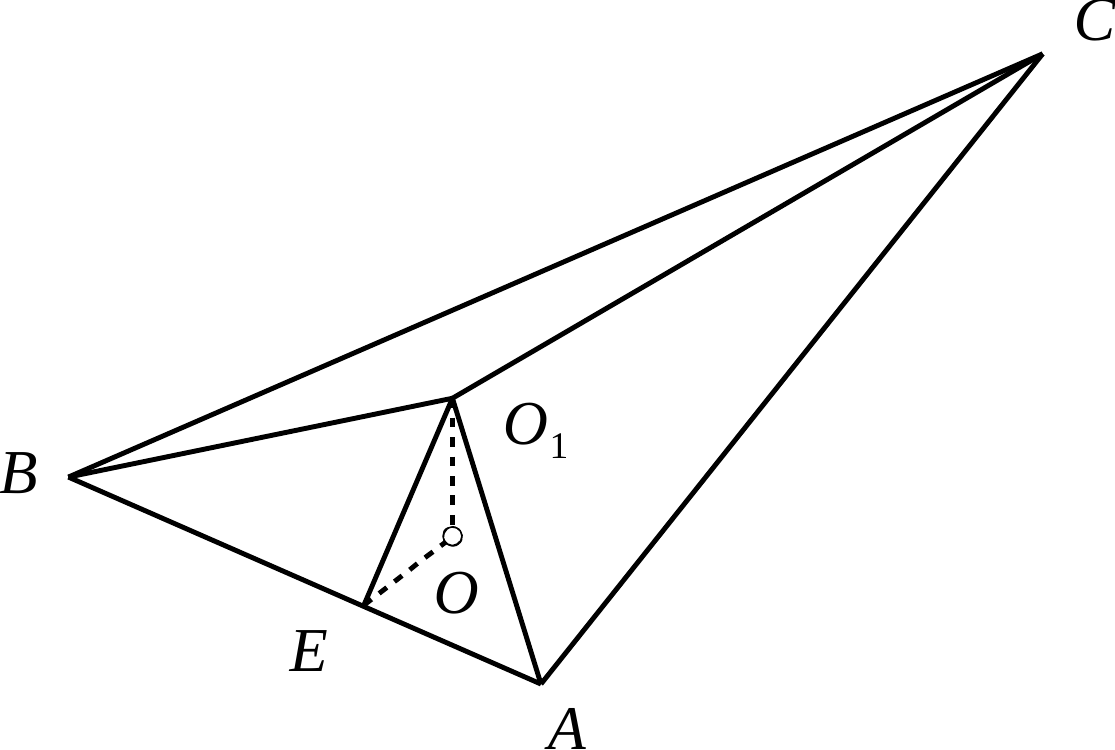}
 \caption{Scheme of the variation $(\delta x).$}
 \label{fig1}
\end{figure}

Let us fix a matrix of rank one: $(\xi^{\alpha}\eta^i)$, where 
$\xi$ and $\eta$, are two unit vectors in the space of variables 
$(x)$ and $(t)$ respectively. Let us build the graph of the variation as a many-dimensional polyhedral  "house of cards"  having the form of thickened pyramid that situated in a neighbourhood of a given point $\tau = (\tau^1,\tau^2, \dots \tau^n)$ of the plane 
$(t)$ (see Fig. \ref{fig1}).

Its projection $\Xi$ on the space of coordinates $(t)$, that corresponds to the support of the variation to be build, is the thickened isosceles triangle $(ABC), \; |AC|=|BC|$. The word "thickened" means that each point on the scheme is the cube with 
the center at this point and with
the edge $\sigma^{2(n-2)}$, and each segment is the parallelepiped with exactly one parallel group of edges which equal to the length of the segment and others equal $\sigma^{2(n-2)}$. The length of the face $|AB|=2\sqrt{\sigma}$. Its center is in the point $E=\tau$. Take the point $O$ lying on the internal normal to $(AB)$ passing through the point $E$ and having coordinate $\{ \tau^i+\frac{\sqrt{7}}{3}\sigma \eta \}$. Take at the point $O$ the perpendicular to the plane $(t)$ at the direction $\xi$ and mark on it the point $O_1$ at the distant 
$\sigma$ from $O$.

The upper point $O_1$ is joined with the vertices of the base by faces $(O_1A)$, $(O_1B)$ and $(O_1C)$. The front face of the base (depicted by the line $(AB)$) is situated on the distant 
$|OE|=\sigma \eta \sqrt{7}/3$ (again without regard to the thickness) from the point $O$. In what follows, all measures (volume, length and so on) will be given without regard for the constant thickness 
$\sigma^{2(n-2)}$.

The external normal vector to $(AB)$ in the plane $(t)$ is $\eta$. We will call the upper front face $\goth M = (AO^1B)$ by the main one. The remaining faces ${\goth M}'$ are called by the minor ones. The length of $|AB|=\sigma^{1/2}$. So, the slope of the main face 
$\goth M = (AO_1B)$ is $(\sqrt{7}/3)\xi^{\alpha}\eta^i$. The face 
$(O_1C)$ is projected into $(OC)$; the face $(O_1A)$ is projected into $(OA)$; the face $(O_1B)$ is projected into $(OB)$.
The length of $|OC|=\sigma^{1/4}$. All the picture has the mirror symmetry relative to the plane $(CO_1EO)$. The main face $\goth M$, having the constant slope, tends to zero in $C$-metric; its equation is $x=(\xi^{\alpha}\eta^i)(\sqrt{7}/3)(t^i-\tau^i)$. Now we have one more argument in favour of using rank one matrices. If the rank would be greater than one, it would be impossible to build the polyhedral continuous variation, since walls with greater rank would intersect the base by planes with co-dimension greater than $1$.

It is easy to see that minor faces $\goth M'$ tend to the horizontal direction uniformly in $C^1$-metric, because its height $|O_1O|=\sigma$ tends to zero much faster than 
$|O_1A|=|O_1B| \sim \sigma^{1/2}$. 

The variation obtained over $\Xi$ will be denoted by 
$(\delta x)(t,\sigma)$, or simply $(\delta x)$.

\bigskip 

\bigskip

\subsection{Variation of the functional}

\bigskip

Consider the solution $\hat x(t)$ with the fixed solution $q(t)$ of system $(\ref{5})$ and with the variation $(\delta x)(t,\sigma)$.
To prove the necessary condition of optimality it is usual to find the derivative of ${\cal F}$ in the direction $(\delta x)(t,\sigma)$
at $\sigma = 0$. In our case there exists an additional difficulty.
The first derivative of ${\cal F}$ is zero 
${\cal F}'|_{\sigma}(+0)=0$ and the second one equals to infinity
${\cal F}''|_{\sigma}(+0)=\infty$. We are forced to find the asymptotic behaviour of ${\cal F}(\sigma)$ as $\sigma \to 0$. 
Let us denote by $\bar{\eta}=\eta \sqrt{7}/3$.
The increment of the integral on the variation in question has the form of the sum of two integrals: on the main face $\goth M$ and on the rest of the polyhedron $\goth M'$                          

$$
\begin{array}{c}
{\cal F}(\hat{x}+\delta x) - {\cal F}(\hat x) =
\int_{\goth M}\left[ f\left( t,\hat x+(\delta x), \hat {\dot x} +\xi^{\alpha}\bar{\eta}_i\right)-f\left( t,\hat x,\frac{D\hat x}{Dt}\right)\right] dt^{\I} \\
+ \int_{\goth M'}\left[ f\left( t, \hat x+(\delta x), \hat {\dot x}+\frac{D(\delta x)}{Dt}\right)  
- f\left( t,\hat x,\frac{D\hat x}{Dt}\right)\right] dt^{\I}
\end{array}
$$

From the condition of the strong minimum follows that
the increment of the integral by $\sigma$ in $\sigma =0$ 
should be non-negative. Consider first the integral over 
$\goth M$. 

\bigskip

Denote by $v(Y)$ the volume of the face $Y$. The volume of the face 
$\goth M$ is equal 
$v(AO_1B)=|O_1E|\cdot |AE|=4/3\sigma^{3/2}$.
The faces $(AO_1) \cup (O_1B)$ will be denoted by ${\goth P}$.
The volume of the faces ${\goth P}$ is equivalent to 
$|{\goth P}| \sim 2\sigma^{1/2}$.

\bigskip

$$
\int_{\goth M}\left[ f\left( t,\hat x+(\delta x), \hat{\dot x}+\xi^{\alpha}\bar{\eta}_i\right) - \hat f\right] dt^{\I}
$$
Using the mean value theorem we obtain
\beq
\begin{array}{c}\Delta {\cal F}_{\sigma}(\goth M)=
\left[f\left(\bar t,\hat{x}(\bar{t})+(\delta x)(\bar t,\sigma),\hat{\dot{x}}(t)+
(\xi^{\alpha}\bar{\eta}_i)\right ) - \hat f(\bar t,\bar x)
\right] v(\goth M) = \\
4/3(\sigma)^{3/2}\left[ f\left( 
\bar t,\hat{x}(\bar{t})+\delta x(\bar t,\sigma),
\dot{\hat{x}}(\bar{t})+
(\xi^{\alpha}\bar{\eta}_i)\right ) - \hat f\right]
\end{array}
\label{14}
\eeq
\noindent where $(\bar t \in {\goth M}$ and $\bar x=\hat x(\bar t)$  is an intermediate value of (t, x), obtained from the mean value theorem.

Let us differentiate the integral over $\goth M'$.

$$
\lim_{\sigma \to +0}\frac{1}{\sigma}\left[
\int_{\goth M'}f\left( t,\hat{x}+(\delta x), \hat{\dot{x}}+
\frac{D(\delta x)(t, \sigma)}{Dt}\right ) - \hat f \right] dt^{\I}.
$$
The differentiation under the sign of the integral gives

$$
\int_{\goth M'}
\frac{\partial \hat f}{\partial x^{\alpha}}\left(
\frac{\partial (\delta x)^{\alpha}}{\partial \sigma}\right) +
\frac{\partial \hat f}{\partial 
\left( \frac{\partial x^{\alpha}}{\partial t^i}\right)} \frac{\partial }{\partial t^i}
\left( \frac{\partial (\delta x)^{\alpha}}{\partial \sigma}\right)
dt^{\I}
$$

Integrate by part the second part of the last integral. Instead of it we have

$$
\int_{\partial {\goth M'}}
\frac{\partial \hat f}{\partial \left( 
\frac{\partial (\delta x)^{\alpha}}{\partial t^i}\right)}
\left( \frac{\partial (\delta x)^{\alpha}}{\partial \sigma}\right) dt^{I_{n-1}} - \int_{\goth M'}\frac{\partial }{\partial t^i}\left( 
\frac{\partial \hat f}{\partial \left( \frac{\partial (\delta x)^{\alpha}}{\partial t^i}\right)}\right)\left( 
 \frac{\partial (\delta x)^{\alpha}}{\partial \sigma}\right)dt^I.
$$
Under the sign of the integral on the domain ${\goth M'}$ it appears the left hand side of the Euler equation that equals zero. So, it remains only the integrals over the boundary of ${\goth M'}$. Integrals along faces $(AC)$ and $(BC)$ are zero, since
$(\delta x)$ on it is zero. It remains the integral along the moving face ${\goth P}$.

We have to rearrange

\beq
-\int_{\goth P}
\frac{\partial \hat f}{\partial \left( 
\frac{\partial (\delta x)^{\alpha}}{\partial t^i}\right)}
\left( \frac{\partial (\delta x)^{\alpha}}{\partial \sigma}\right) dt^{I_{n-1}}.  
\label{15}
\eeq
The sign "minus" is due to the fact that the orientation of the face ${\goth P}$ which is induced by the face $\goth M'$ is opposite to the orientation that was induced  by the face $\goth M$.
As it is usual for variational equations the value of
$\frac{\partial (\delta x)}{\partial \sigma}$ is calculated relative to the fixed values of $t$. But here we need to differentiate
along face $\goth M'$ with moving boundary. So, we need to use the full derivative $\frac{d(\delta x)(t,\sigma)}{d\sigma}$. It should be taken into account that $\frac{dt}{d\sigma}=\sqrt{7}/3$. The differentiation gives
$$
\frac{d(\delta x)(t,\sigma)}{d\sigma}=
\frac{\partial (\delta x(t,\sigma))}{\partial \sigma}+
\frac{D(\delta x(t,\sigma))}{Dt}\sqrt{7}/3.
$$
It will be recalled that on the moving face $(\delta x)(t,\sigma) =
(\xi^{\alpha}\bar{\eta}^i)(t^i-\tau^i)$. Hence, 
$$
\frac{\partial (\delta x)(t,\sigma)}{\partial \sigma}=
\frac{d(\delta x)(t,\sigma)}{d\sigma}-
(\xi^{\alpha}\bar{\eta}_i).
$$
We rearrange $(\ref{15})$ by using the last formula. We obtain

$$
-\int_{\goth P}\frac{\partial \hat{f}}{\partial \left( 
\frac{\partial (\delta x)^{\alpha}}{\partial t^i}\right)}
(\xi^{\alpha}\bar{\eta}^i) dt^{I_{n-1}} + 
\int_{\goth P}\frac{\partial \hat{f}}{\partial \left( 
\frac{\partial (\delta x)^{\alpha}}{\partial t^i}\right)}
\frac{d(\delta x)^{\alpha}}{dt}dt^{I_{n-1}}. 
$$

\bigskip
To restore the increment we have to integrate the expression by
$\sigma$. We again use the mean value theorem to obtain the following expression for the integral.
\beq
\Delta {\cal F}_{\sigma}(\goth M')=
4/3(\sigma)^{3/2}\left[ -
\frac{\partial \hat{f}}{\partial \left( 
\frac{\partial (\delta x)^{\alpha}}{\partial t^i}\right)}(\tilde{t},
\tilde{x})(\xi^{\alpha}\bar{\eta}_i)  + 
\frac{\partial \hat{f}}{\partial \left( 
\frac{\partial (\delta x)^{\alpha}}{\partial t^i}\right)}
\frac{d(\delta x)^{\alpha}}{dt}(\tilde{t}, \tilde{x})
\right]
\label{16}
\eeq
\noindent Here $(\tilde{t}, \tilde{x}) \in {\goth P}$ is the point chosen by the mean value theorem.

\bigskip
Both points $(\tilde{t}, \tilde{x})$ and $(\bar{t}, \bar{x})$ tend to the central point $(\tau , \hat{x}(\tau))$ in the process of shrinking of the polyhedron as $\sigma \to 0$. Both coefficients $4/3(\sigma)^{3/2}$ before the expressions of
$\Delta {\cal F}_{\sigma}$ for both faces $\goth M$ and $\goth M'$ are the same. Consequently, we have the following asymptotic formula
for any matrix $(\xi^{\alpha}\bar{\eta}_i)$ of the rank 1:
$$
f\left( t,\hat{x},\frac{D\hat{x}}{Dt}+(\xi^{\alpha}\bar{\eta}_i ) \right) - f\left( t,\hat{x},\frac{D\hat{x}}{Dt}\right)-
\left( \frac{\partial \hat{f}}{\partial (\frac{\partial x^{\alpha}}
{\partial t^i})} \right) (\xi^{\alpha}\bar{\eta}_i) +
\left( \frac{\partial \hat{f}}{\partial (\frac{\partial x^{\alpha}}
{\partial t^i})} \right)\left( \frac{D\hat{x}}{Dt}\right) \leq 0.
$$

In other words, in terms of the Pontryagin's function     

$$
\max_{(\xi,\eta)\in {\cal U}}{H}\left( t,\hat x,\hat{\dot{x}}+
\xi^{\alpha}\eta_i) \right) = {H}\left( t,\hat x,\hat{\dot{x}}\right).
$$

Thus, the theorem is proved for any internal points of the domain 
$V$. At points of boundary it can be proved by the passage to the limit in corresponding inequalities.

$\Box$

\begin{example}

By way of example, consider a problem associated with elasticity theory

To minimize the functional

\beq
\int_{\goth N}[a(z_1^2+z_4^2)+b(z_2^2+z_3^2)+2c\det z]
dt^1\wedge dt^2
\label{17}
\eeq

Here $z_1=\frac{\partial x^1}{\partial t^1},\,
z_2=\frac{\partial x^1}{\partial t^2},\,
z_3=\frac{\partial x^2}{\partial t^1},\,
z_4=\frac{\partial x^2}{\partial t^2}.
$  The summand, containing $\det z = z_1z_4-z_2z_3$, defines
the degree of contraction-expansion of the material; the coefficient $2c$ is called by solid elasticity modulus. The coefficients $a,b$ binds to the constants of Lame, that is to the tensor of elasticity module. The slopes of the surface (variables $z_i$) serve as control variables.

The matrix of the quadratic form, standing under the sign of the integral, is

$$
\left(
\begin{array}{cccc}
a&  0&   0&   c\\
0&  b&  -c&   0\\
0& -c &  b&   0\\
c&  0&   0&   a\\
\end{array}
\right)
$$
Eigenvalues of the matrix are 
$\lambda_1=a-c,\, \lambda_2=a+c,\, \lambda_3=b-c,\, \lambda_4 =b+c.$ Say $a>b$. For $a>c>b$ one of the
eigenvalues is negative. The quadratic form is non-convex. The maximum of the Pontryagin's function
$$
H=-\left\{
a\left(
\frac{\partial x^1}{\partial t^1}\right)^2+
a\left(
\frac{\partial x^1}{\partial t^4}\right)^2+
b\left(
\frac{\partial x^2}{\partial t^1}\right)^2+
b\left(
\frac{\partial x^1}{\partial t^2}\right)^2+
2c\left(
\frac{\partial x^1}{\partial t^1}\frac{\partial x^2}{\partial t^4}-
\frac{\partial x^1}{\partial t^2}\frac{\partial x^2}{\partial t^1}
\right) \right\} + q^i_{\alpha}
\frac{\partial x^{\alpha}}{\partial t^i}
$$
\noindent (that would be in line with the naive generalization of the Pontryagin's maximum principle) 
does not attained on any extremal. Nevertheless, the restriction of this function to the level-surface of the rank one matrices reduces $H$ to the positive definite quadratic form.

Indeed, the first variation on extremals is zero. Variations of control $h=(h_1,\, h_2,\, h_3,\, h_4)$, which correspond to directions with the rank one matrices, are equivalent to the degenerate matrices $h$. So, the main quadratic part of expansion relative to $h$ of the summand $\det z$ on such variations (that is $h_1h_4-h_2h_3$) equals zero. It remains the positive definite quadratic form which ensure our maximum principle. To test sufficient conditions it should be appealed to the theory of fields of extremals that was developed in  $\cite{MZ}$. 

\end{example}
 
\begin{example}
 
Consider the problem of minimization of the functional

\beq
\int_{\goth N}[(z_1)^3+(z_2)^3)]dt^1\wedge dt^2.
\label{18}
\eeq

Here $z_1=\frac{\partial x^1}{\partial t^1},\,
z_2=\frac{\partial x^2}{\partial t^2}$; 
$V = \{ t_1^2+t_2^2 \le 1\}$. The boundary conditions are 
$x_1|_V = \cos \varphi ,\quad  x_2|_V = \sin \varphi$.

The Euler equations have the form

\beq
\left\{ 
\begin{array}{c}
\frac{\partial}{\partial t_1} 
\left( 3\left(\frac{\partial x_1}{\partial t_1}\right)^2\right) = 0, \\
\frac{\partial}{\partial t_2} 
\left( 3\left(\frac{\partial x_2}{\partial t_2}\right)^2\right) = 0 
\end{array}
\right.
\label{19}
\eeq

It is ease to see that the unique solution to Euler equations satisfying the boundary conditions is
$x_1=t_1; \; x_2=t_2$. The second variation at the extremal on matrices of the rank 1 equals
$6\xi_1^2\eta_1^2 + 6\xi_2^2\eta_2^2$. It is strictly positive,
so the condition of Hadamard-Legendre for the weak minimum is fulfilled. But the Pontryagin's function on matrices of the rank 1
 equals $-(\xi_1^3\eta_1^3+\xi_2^3\eta_2^3) + 
6(\xi_1\eta_1+\xi_2\eta_2)$. It reaches only local maximum on the extremal $\{\xi_1\eta_1=1, \; \xi_2\eta_2=1\}$. The global maximum equals $+\infty$. By theorem 2 we conclude that the minimum on the extremal is not strong.    

\end{example}

\end{document}